\documentclass[epsfig,12pt]{amsart}
\usepackage{amsmath, amsthm, amsfonts, amssymb, color, graphicx}
\usepackage{mathrsfs}
\newtheorem{thm}{\sc Theorem}
\newtheorem{lem}{\sc Lemma}
\newcommand{\scr}[1]{\mathscr #1}

\def\N{\mathbb N} 
\def\Re{{\rm Re}\;} 
\def\Im{{\rm Im}\;}
\def\D{\mathbb{D}}
\def\J{\scr J} 

\def\A{\scr{A}} 
\def\W{\mathbf{W}}
\def\B{\scr B}
\def\BB{\mathbf{B}}
\def\OM{\mathbf{\Omega}} 
\def\C{\mathbb{C}}     
\def\ds{\displaystyle}
\def\kzu1#1{\buildrel #1 \over \longrightarrow}

\begin{document}

\begin{center}{\large\sc On the Dynamics of Rational Maps with Two Free Critical Points}\\
\sc HyeGyong Jang\footnote{Faculty of Mathematics, University of Science {Pyongyang} D.P.R. of Korea.}\footnote{This paper was written during a visit of the CAS supported by the TWA-UNESCO Associateship Scheme.} and Norbert Steinmetz\footnote{Fakult\"at f\"ur Mathematik, TU  Dortmund, Germany.\\ http://www.mathematik.uni-dortmund.de/steinmetz/ \\ e-mail: stein$@$math.tu-dortmund.de}\\
\today\end{center}

\bigskip
\begin{quote}\footnotesize {\sc Abstract.} In this paper we discuss the
dynamical structure of the rational family $(f_t)$ given by
$$f_t(z)=tz^{m}\Big(\frac{1-z}{1+z}\Big)^{n}\quad(m\ge 2,~t\ne 0).$$
Each map $f_t$  has two super-attracting immediate basins and two free critical points. We prove that for $0<|t|\le 1$ and $|t|\ge 1$, either of these basins is completely invariant and at least
one of the free critical points is inactive.
Based on this separation we draw a detailed picture the structure of the dynamical and the parameter plane.\end{quote}

\section{Introduction}
\noindent Non-trivial rational families $(f_t)$ normally contain specific maps of different character with most interesting and unexpected Julia sets:

\begin{itemize}
\item[-] totally disconnected Julia sets (Cantor sets) occur in any family $z\mapsto z^d+t$;
\item[-] Julia sets consisting of uncountably many (a Cantor set of) quasi-circles occur in the McMullen family $z\mapsto z^m+t/z^n$, which was introduced in \cite{McMullen}. The number of papers on various features of this family is legion; \cite{Devaney} marks the preliminary end of a long list of papers.
\item[-] Julia sets that are Sierpi\'nski curves (Tan Lei and J.\ Milnor~\cite{TLMilnor} were the first to construct examples with this property) again in the McMullen family \cite{Ste1},
the  Morosawa-Pilgrim family $z\mapsto t\big(1+\frac{(4/27)z^3}{1-z}\big)$ \cite{DFGJ,Ste3}, and the family $t\mapsto-\frac t4\frac{(z^2-2)^2}{z^2-1}$ \cite{JS}.
\item[-] In any reasonable family, copies of the Mandelbrot sets of the families $z\mapsto z^d+t$ are dense in the bifurcation locus -- the Mandelbrot set is universal \cite{McMullen00}.
\end{itemize}

Each of these families has just one {\it free} critical point (or several free critical points which have the same dynamical behaviour, this happens, for example, in the McMullen family; the quasi-conjugate family $F_t(z)=z^m(1+t/z)^d$ has just one free critical point). In contrast to this the rational maps
\begin{equation}\label{TheFamily}
f_t(z)=tz^m\Big(\frac{1-z}{1+z}\Big)^n\quad(m\ge 2,~n\in\N,~d=m+n,~t\ne 0)
\end{equation}
in the family under consideration have two free critical points. In this paper we will give a complete description of the parameter plane and the various dynamical planes.
For basic notations and results the reader is referred to the texts \cite{Bea,Cal,Mil,McMullen,Ste}.
\section{Notation}

\noindent The rational map (\ref{TheFamily}) has
\begin{itemize}
\item[-]
two super-attracting fixed points $0$ and $\infty$ with corresponding basins $\A_t$ and $\B_t$, respectively. Then either $\A_t$, say, is completely invariant or else has a single pre-image $\A_t^*$ that is mapped in a $(n:1)$-manner onto $\A_t$, which will be written as
$$\A_t^*\kzu1{n:1}\A_t;$$
\item[-] two free critical points
$$\qquad\alpha=-\frac nm+\sqrt{1+\Big(\frac nm\Big)^2}\quad{\rm and}\quad \beta=-\frac nm-\sqrt{1+\Big(\frac nm\Big)^2}$$
and critical values
$$v^\alpha_t=f_t(\alpha)=tv^\alpha_1\quad{\rm and}\quad v^\beta_t=f_t(\beta)=t v^\beta_1;$$
\item[-] two {\it escape loci} $\OM^\alpha$ and $\OM^\beta$, with  $t\in\OM^\alpha$ and $t\in\OM^\beta$ if and only if $f^k_t(\alpha)\to 0$ and $f^k_t(\beta)\to \infty$, respectively, as $k\to\infty$;
\item[-] two {\it residual sets} $\OM_{\rm res}^\alpha$ and $\OM_{\rm res}^\beta$, with $t\in\OM_{\rm res}^\alpha$ and $t\in\OM_{\rm res}^\beta$ if and only if $v^\beta_t\in\A_t$
and $v^\alpha_t\in\B_t$, respectively.
\end{itemize}
The notation of the residual sets indicates that $\OM_{\rm res}^\alpha$ is related to $\OM^\alpha$ rather than $\OM^\beta$.
The open sets $\OM^\alpha$ and $\OM^\beta$ are in a natural way sub-divided into
\begin{itemize}
\item[-] $\OM_0^\alpha$ resp. $\OM_0^\beta$: $v^\alpha_t\in\A_t$ resp. $v^\beta_t\in\B_t$, and
\item[-] $\OM_k^\alpha$ resp. $\OM_k^\beta$: $f^k_t(v^\alpha_t)\in\A_t$, but $f^{k-1}_t(v^\alpha_t)\notin\A_t$ resp.\\
\phantom{$\OM_k^\alpha$ resp. $\OM_k^\beta$:} $f^k_t(v^\beta_t)\in\B_t$, but $f^{k-1}_t(v^\beta_t)\notin\B_t$ ($k\ge 1)$.
\end{itemize}

Hitherto, $f_t$ is hyperbolic and the Fatou set of $f_t$ consists of the basins $\A_t$ and $\B_t$, and their pre-images, if any. However, there may and will be
also other hyperbolic components. By $\W_k^\alpha$ and $\W^\beta_k$ we denote the open sets such that $\alpha$ and $\beta$ belongs to some
(super-)attracting cycle of Fatou domains $U_1,\ldots,U_k$, respectively, not containing $0$ and $\infty$.

The {\it bifurcation} locus $\BB$ of the family $(f_t)_{0<|t|<\infty}$ is the set of $t$ such that the Julia set $\J_t$ does not move continuously over any neighbourhood of $t$, see McMullen~\cite{McMullen}. In order that $t\in\BB$ it is necessary and sufficient that at least one of the free critical points
is {\it active}. Thus $\BB=\BB^\alpha\cup\BB^\beta$, where $t\in\BB^\alpha$ resp. $t\in\BB^\beta$ means that $\alpha$ resp. $\beta$ is active. It is {\it a priori} not excluded that $\BB^\alpha$ and $\BB^\beta$ overlap. Although there is just one parameter plane, each point of this plane carries at least two pieces of information, so one could also speak of the $v^\alpha_t$- and $v^\beta_t$-plane.

We also set
$$Q_0(t)=v^\alpha_t=tv^\alpha_1\quad{\rm and}\quad Q_k(t)=f_t^k(v^\alpha_t)=f_t(Q_{k-1}(t))\quad(k\ge 1)$$
and note that $Q_k$ is a rational function of degree $1+d+\cdots +d^k=\frac{d^{k+1}-1}{d-1}$ with a zero of order $\frac{m^{k+1}-1}{m-1}$ at the origin.

From
$$-1/f_t(-1/z)=f_{(-1)^{d+1}/t}(z)\quad(d=m+n)$$
it follows that $f_t$ is conjugate to $f_{1/t}$ if $d$ is odd, and to $f_{-1/t}$ if $d$ is even, hence $\OM^\alpha=\pm 1/\OM^\beta$, and this is also is true for $\OM^\alpha_{k}$ and $\OM^\beta_{k}$, $\OM^\alpha_{\rm res}$ and $\OM^\beta_{\rm res}$, $\W_k^\alpha$ and $\W_k^\beta$, and $\BB^\alpha$ and $\BB^\beta$. This also indicates that the circle $|t|=1$ plays a distinguished role
with strong impact on what follows.

\begin{lem}\label{TRENNUNG}For every $m\ge 2$, $n\ge 1$ there exists some $r>0$, such that for $0<|t|\le 1$ the disc $\triangle_{r|t|}:|z|<r|t|$ contains $f_t(\overline\triangle_{r|t|}\cup [0,1])$, but does not contain $v^\beta_t$.\end{lem}

\proof We will first consider $f_1$ and show that there exists some disc $\triangle_{r}:|z|<r$ such that $f_1(\overline\triangle_r\cup [0,1])\subset\triangle_r$
holds. This is easy to show if  $n<m$ for $r=\frac 13$:
$$|f_1(z)|\le 3^{-m}2^n<\textstyle\frac13$$
holds if $|z|\le \frac13$ and $m>n\ge 1$, and from
$$0\le f_1(x)\le x^2\frac{1-x}{1+x}\le\textstyle \frac12(5\sqrt 5-11)<\frac1{10}\quad(0\le x\le 1)$$
the assertion follows.\\
We now consider the case $n\ge m$. Then $f_1(\overline\triangle_r)\subset\triangle_r$ holds as long as
$$g(r)=r^{m-1}\Big(\frac{1+r}{1-r}\Big)^n<1,$$
and $f_1$ maps $[0,1]$ into $\triangle_r$ provided
$$v^\alpha_1=\max_{0\le x\le 1}x^m\Big(\frac{1-x}{1+x}\Big)^n<r.$$
Since $g$ is increasing this may be achieved if $g(v^\alpha_1)<1$ holds. To prove this we note that
$\sqrt{1+\tau}-1=\frac\tau{2\sqrt{1+\theta\tau}}$ $(0<\theta<1,$ $\tau=\frac{m^2}{n^2}\le 1)$
implies $\frac m{2\sqrt 2n}<\alpha<\frac m{2n}$, while from $\log\frac{1-x}{1+x}< -2x$ ($0<x<1$) it follows that
$$v^\alpha_1<\Big(\frac m{2n}\Big)^me^{-2\frac{m}{2\sqrt 2}}=\Big(\frac{m}{2e^{\frac1{\sqrt 2}} n}\Big)^m<\Big(\frac m{4n}\Big)^m=\mu^m.$$
Moreover, from
$$\log\frac{1+x}{1-x}=2x\Big(1+\frac 13x^2+\frac15x^4+\cdots\Big)\le 2x\Big(1+\frac{x^2}3\frac 1{1-x^2}\Big)\le\textstyle  2x(1+\frac1{45}),$$
which holds for $x=\big(\frac m{4n}\big)^{m-1}\le \frac 1{4}$, we obtain
$$\left(\frac{1+\mu^m}{1-\mu^m}\right)^n=\left(\frac{1+\frac m4\frac{\mu^{m-1}}{n}}{1-\frac m4\frac{\mu^{m-1}}{n}}\right)^n\le e^{\frac{23}{45}m\mu^{m-1}}<\Big(e^{(\frac m{4n})^{m-1}}\Big)^m.$$
Thus $g(v^\alpha_1)<1$ follows from $\big(\frac m{4n}\big)^{m-1}e^{(\frac m{4n})^{m-1}}\le \frac 1{4}e^{\frac 14}<1$.

With this choice of $r\in(0,1)$ it is clear that $v^\beta_t$ belongs to $\triangle_r$ if $|t|$ is small. For individual $0<|t|\le 1$, $f_t(z)=tf_1(z)$ maps $\overline\triangle_{r|t|}\cup[0,1]$ into $\triangle_{r|t|}$, while $v^\beta_t\notin\triangle_{r|t|}$ follows from $|v^\beta_t|=|t|/v^\alpha_1>|t|>r|t|$.
\hfill$\Box$

\section{The escape loci}
\noindent The purpose of  Lemma~\ref{TRENNUNG} is twofold. First of all it shows that the critical points $\alpha$ and $\beta$ cannot be simultaneously active, and the bifurcation sets $\BB^\alpha$ and $\BB^\beta$ are separated by the unit circle $|t|=1$.
Secondly, the condition $v^\beta_t\notin\triangle_{r|t|}$ ($0<|t|\le 1)$ ensures that in an exhaustion $(D_\kappa)$ of $\A_t$ starting with $D_0=\triangle_{r|t|}$, $D_\kappa$ is simply connected
as long as $\beta\notin D_\kappa$, and $f_t:D_{\kappa}\kzu1{d:1} D_{\kappa-1}$ has degree $d=m+n$. In particular, for $t\in\OM^\alpha_{\rm res}$ there exists some simply connected and forward invariant domain $D_\kappa\subset\A_t$ that contains $v^\beta_t$.\medskip

We note some more simple consequences of Lemma~\ref{TRENNUNG}; our focus is on the critical point $\alpha$ and the ${}^\alpha$-sets.
\begin{itemize}
\item[-] $\{t:0<|t|\le 1\}\subset\OM^\alpha_0$;
\item[-] $\overline{\OM^\alpha_{\rm res}}\subset\D;$
\item[-] $\alpha$ is inactive on $0<|t|\le 1$;
\item[-] $\overline{\bigcup_{k\ge 1}(\OM^\alpha_k\cup\W^\alpha_k)}\subset\{t:1<|t|<T\}$ for some $T=T_{mn}>1$;
\item[-] $\BB^\alpha\subset\{t:1<|t|<T\}$ for some $T=T_{mn}>1$.
\end{itemize}

The consequences for the dynamical planes are as follows.

\begin{thm}\label{OM0capOMinfty}For $t\in\OM^\alpha_0$, the basin $\A_t$ is completely invariant, and any other Fatou component is simply connected. Moreover,
\begin{itemize}
\item[-] for $t\in\OM^\alpha_0\cap\OM^\beta_0$ also $\B_t$ is completely invariant, the Julia set $\J_t=\partial\A_t=\partial\B_t$ is a quasi-circle, and $f_t$ is quasi-conformally conjugate
to $z\mapsto z^d$;
\item[-] for $t\in\OM^\alpha_{\rm res}$, $\A_t$ is infinitely connected and the Fatou set consists of $\A_t$, $\B_t$, and the predecessors of $\B_t$ of any order.
\end{itemize}
\end{thm}

\proof To prove complete invariance of $\A_t$ we first assume $0<|t|\le 1$. Then $\A_t$ contains the interval $[0,1]$ by Lemma~\ref{TRENNUNG}, hence is completely invariant.
If, however, $|t|>1$, then $\B_t$ is completely invariant, and any other Fatou component is simply connected.
Assuming $1\not\in\A_t$ ($t\in\OM^\alpha_0$, $|t|>1$) we obtain either $f_t:\A_t^*\kzu1{n:1}\A_t$ with $n=(n-1)+1$ critical points if $\alpha\in\A_t^*$ or else
$f_t:\A_t\kzu1{m:1}\A_t$ with $m=(m-1)+1$ critical points if $\alpha\in\A_t$, this contradicting simple connectivity of both domains $\A_t$ and $\A_t^*$ by the Riemann-Hurwitz formula.

The first assertion is obvious since $\B_t$ shares the properties of $\A_t$ and $f_t$ is hyperbolic.

The second assertion follows from the Riemann-Hurwitz formula, since $f_t:\A_t\kzu1{d:1} \A_t$
has degree $d$ and $r=(m-1)+(n-1)+1+1=d$ critical points $0$, $1$ (if $n>1$), $\alpha$, and $\beta$.\hfill$\Box$

\begin{center}
\includegraphics[scale=0.133]{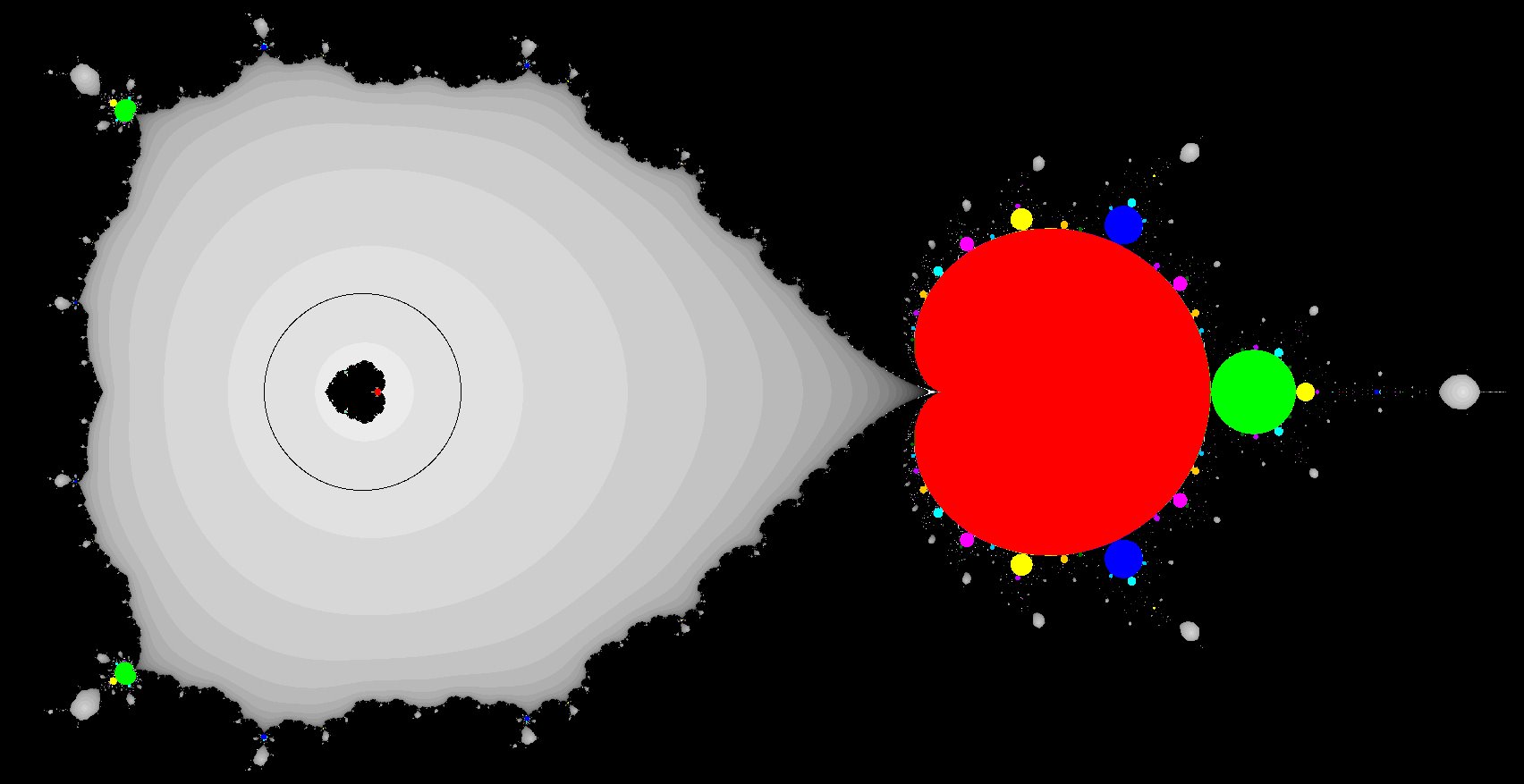}~
\includegraphics[scale=0.1526]{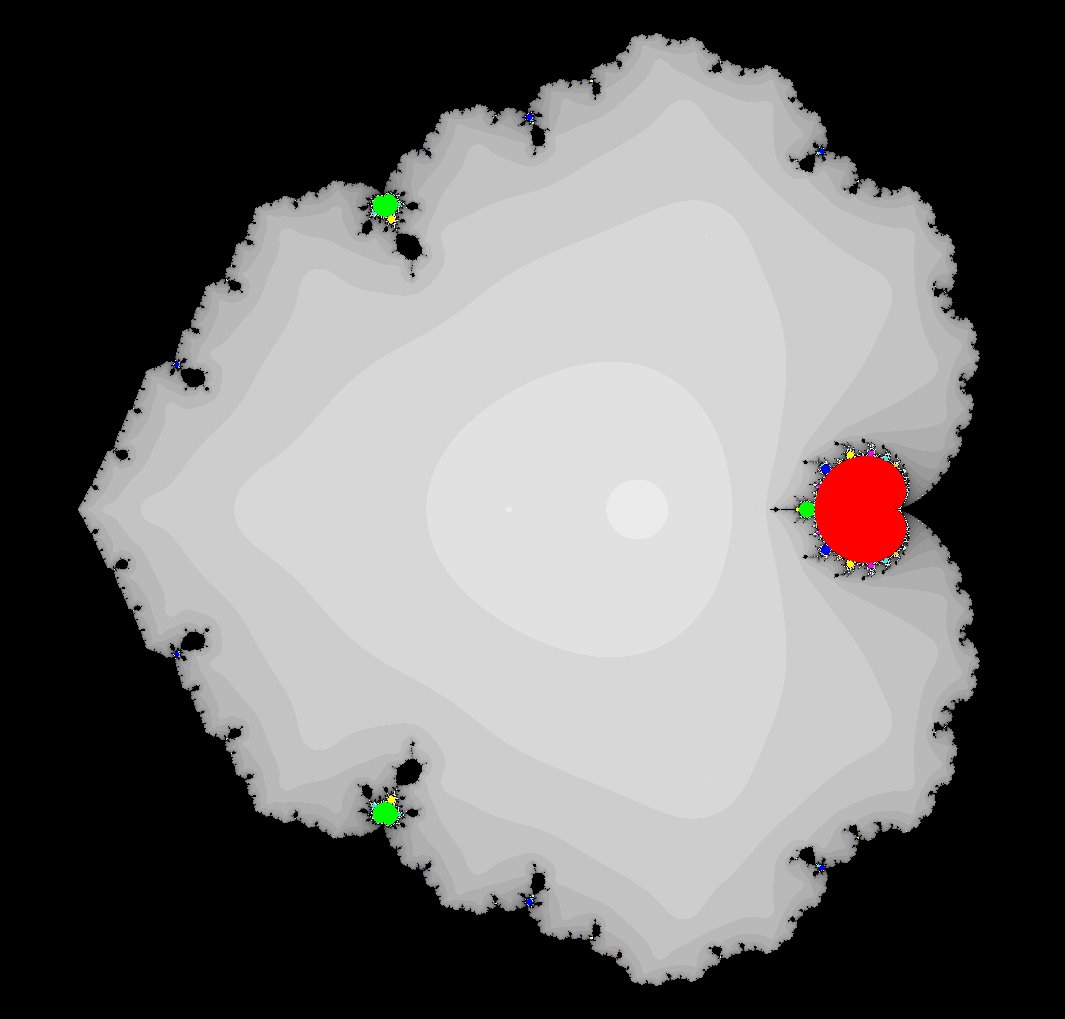}

\medskip\footnotesize {\sc Figure~1}-left: The $\alpha$-parameter plane for $\ds f_t(z)=tz^2\frac{1-z}{1+z}$ displaying the unit circle, $\OM^\alpha$ in gray, $\OM^\alpha_{\rm res}$ and $\OM^\beta_{\rm res}$ in black (in- and outside the unit circle), and $\W^\alpha$ (coloured).
{\sc Figure~1}{-right:} A neighbourhood of the origin displaying $\OM^\alpha_{\rm res}$ surrounded by points of $\OM^\alpha_0$ (black), $\OM^\beta_k$ ($k\ge 1$, black, small), and  $\W^\beta$ (coloured). \end{center}

\begin{thm}\label{SimpleConn}$\OM^\alpha_0\cup\{0\}$, $\OM^\alpha_{\rm res}\cup\{0\}$, and the connected components of $\OM^\alpha_k$ $(k\ge 1)$ are simply connected domains.
Riemann maps onto $\D$  are given by any branch of $\sqrt[m]{E_0(t)}$,  $\sqrt[m]{E_{\rm res}(t)}$, and $\sqrt[n]{E_{k}(t)}$, respectively.\end{thm}

For the proof we need two auxiliary results on the maps

\begin{equation}\label{FunctionsE}\begin{array}{rcll}
E_0(t)&=&t(\Phi_t(v^\alpha_t))^{m-1}&(t\in\OM^\alpha_0),\cr
E_{\rm res}(t)&=&t(\Phi_t(v^\beta_t))^{m-1}&(t\in\OM^\alpha_{\rm res}), {\rm ~and}\cr
E_k(t)&=&t^{\frac 1{m-1}}\Phi_t(f^k(v^\alpha_t))&(t\in\OM^\alpha_k,~k\ge 1),
\end{array}\end{equation}
where $\Phi_t$ denotes the B\"ottcher function to the fixed point $z=0$.
In the first step (Lemma~\ref{Lemma2}) of the proof of Theorem~\ref{SimpleConn} we will show that the functions (\ref{FunctionsE}) provide proper maps on $\D\setminus\{0\}$ and $\D$, respectively, which are only ramified over the origin. In the second step (Lemma~\ref{Lemma3}) this will be used to show that the corresponding domains (with $0$ included, if necessary) are simply connected.

The solution to B\"ottcher's functional equation
\begin{equation}\label{BFEQ}\Phi_t(f_t(z))=t\Phi_t(z)^m\quad(\Phi_t(z)\sim z {\rm~as~}z\to 0)\end{equation}
is locally given by
$$\Phi_t(z)=\lim_{k\to\infty}\sqrt[m^k]{f_t^k(z)/t^{1+m+\cdots +m^{k-1}}}=t^{-\frac 1{m-1}}\lim_{k\to\infty}\sqrt[m^k]{f_t^k(z)};$$
it conjugates $f_t$ to $\zeta\mapsto \zeta^m$. This conjugation holds throughout $\A_t$ in the third case, when $\Phi_t$ maps $\A_t$ conformally onto the disc $|z|<|t|^{-\frac1{m-1}}$; the maps $E_k$ are analytic and well-defined on the components of $\OM^\alpha_k$, $k\ge 1$.

In the first case the conjugation holds on some simply connected neighbourhood of $z=0$ that contains $z=0$ and $z=v^\alpha_t$, but does not contain $z=1$.
The analytic continuation of $\Phi_t$ causes singularities
at $z=1$ and its preimages under $f_t^k$, nevertheless $|\Phi_t(z)|$ is well-defined on $\A_t$ and $|\Phi_t(z)|\to |t|^{-\frac 1{m-1}}$ as $z\to\partial\A_t$ holds anyway.
Thus  $E_0(t)=t\Phi_t(v^\alpha_t)^{m-1}$ is holomorphic on $\OM^\alpha_0$ and zero-free, with $E_0(t)\sim t(v^\alpha_t)^{m-1}=f_1(\alpha)^{m-1}t^m$ as $t\to 0$.

In the second case we construct an exhaustion $(D_\kappa)$ of $\A_t$ such that $f_t:D_\kappa\kzu1{d:1} D_{\kappa-1}$ has degree $d$ and $D_\kappa$ is simply connected for $\kappa\le \kappa_0$ with $v^\beta_t\in D_{\kappa_0}$ and $\beta\in  D_{\kappa_0+1}\setminus D_{\kappa_0}$. This is possible by Lemma~\ref{TRENNUNG}, and the procedure applied to $t^{-\frac 1{m-1}}\Phi_t(v^\alpha_t)$ on $\OM^\alpha_{0}$ also applies to
$t^{-\frac 1{m-1}}\Phi_t(v^\beta_t)$ on $\OM^\alpha_{\rm res}$.

\begin{lem}\label{Lemma2}The functions in (\ref{FunctionsE}) are well-defined and provide proper maps from $\OM^\alpha_0\cup\{0\}$, $\OM^\alpha_{\rm res}\cup\{0\}$, and the connected components of $\OM^\alpha_k$ with $k\ge 1$, respectively, onto the unit disc $\D$.
\end{lem}

\proof To prove that $|E_0(t)|\to 1$ as $t\in\OM^\alpha_0$ tends to $\partial\OM^\alpha_0\setminus\{0\}$ we choose any disc $\triangle:|z|<r$ that is invariant under $f_t$ for every $t\in\OM^\alpha_0$. This is possible since $\OM^\alpha_0$ is contained in some disc $|t|<T$, hence we may choose $r<1$ such that $Tr^{m-1}\big(\frac{1+r}{1-r}\big)^n=1$ holds. By $k=k(t)$ we denote the largest integer such that $f_t^k(v^\alpha_t)\not\in\triangle$. Then $k(t)\to\infty$ as $t\to\partial\OM^\alpha_0\setminus\{0\}$,  and $|f_t^{k(t)}(v^\alpha_t)|\ge r$ implies
$$\liminf_{t\to\OM^\alpha_0\setminus\{0\}} |\Phi_t(v^\alpha_t)|\ge\lim_{t\to\OM^\alpha_0\setminus\{0\}} |t|^{-\frac{1}{m-1}}\sqrt[m^{k(t)}]{r}=|t|^{-\frac{1}{m-1}},$$
while $|\Phi_t(z)|<|t|^{-\frac{1}{m-1}}$ is always true. Thus $E_0$ maps each connected component of $\OM^\alpha_0$ properly onto  $\D\setminus\{0\}$.
It follows that the origin is removable for (a zero of) $E_0$, and $\OM^\alpha_0\cup\{0\}$ is a domain which is mapped by $E_0$ properly with degree $m$ onto the unit disc $\D$.

If $t\in\OM^\alpha_k$ for some $k\ge 1$, then again $|E_k(t)|$ tends to $1$ as $t\to\partial\Omega$, where $\Omega$ is any component of $\OM^\alpha_k$. Thus $E_k$ is a proper map of $\Omega$ onto $\D$. We will prove that $E_k$ is ramified only over zero even for $k\ge 0$, that is $E_k'(t)=0$ implies $E_k(t)=0$. This is a well-known procedure, the idea of which is due to Roesch~\cite{Roesch}, and outlined in detail for the Morosawa-Pilgrim family $z\mapsto t\big(1+\frac{(4/27)z^3}{1-z}\big)$ in \cite{Ste3}, Lemma~2.

We take any $t_0\in\OM^\alpha_k$ and choose
$\varepsilon>0$ such that for $t$ sufficiently close to $t_0$, the closed disc $\triangle_{3\epsilon}:|w-v^\alpha_{t_0}|\le 3\varepsilon$ belongs to the
Fatou component $D_{t_0}$ of $f_{t_0}$ containing $v^\alpha_{t_0}$ ($D_{t_0}$ is a predecessor of $\A_{t_0}$ of order $\ell\ge 0$). Furthermore let
$\eta_t:\widehat\C\longrightarrow \widehat\C$ be any
diffeomorphism such that $\eta_t(w)$ depends analytically on $t$,
$\eta_t(w)=w$ holds on $|w-v^\alpha_{t_0}|\ge 3\varepsilon$ and $\eta_t(w)=w+(v^\alpha_t-v^\alpha_{t_0})$ on
$|w-v^\alpha_{t_0}|<\varepsilon.$ Then
$g_t=\eta_t\circ f_{t_0}:\widehat\C\longrightarrow\widehat\C$
is a quasi-regular map which equals $f_{t_0}$ on $\widehat\C\setminus f^{-1}_{t_0}(\triangle_{3\epsilon})$, and is analytic on $\widehat\C\setminus f_{t_0}^{-1}(A)$ with $A=\{w:\varepsilon\leq|w-v^\alpha_{t_0}|\leq3\varepsilon\}$. To apply Shishikura's qc-lemma \cite{Shi} we need to know that $g_t$ is uniformly $K$-quasi-regular, that is, all iterates $g_t^p$ are $K$-quasi-regular
with one and the same $K$.
This is obviously true if the sets $f_{t_0}^{-p}(A)$ ($p=1,2,\ldots$) are visited at most once by any iterate of $g_t$. This is trivial if $k\ge 1$: the sets $f_{t_0}^{-p}(A)$
belong to different Fatou components, namely  predecessors of $D_{t_0}$ of order $p$. If $k=0$ the argument is different. Let $\triangle_0:|z|<\delta$ be such that $f_{t_0}(\overline{\triangle}_0)\subset\triangle_0$ and set $\triangle_\nu=f_{t_0}^{-1}(\triangle_{\nu-1})$. Then choosing $\epsilon>0$ sufficiently small we have
$A\subset\triangle_{\ell}\setminus\overline{\triangle_{\ell-1}}$ for some $\ell$ and $f_{t_0}^{-p}(A)\subset\triangle_{\ell+p}\setminus\overline{\triangle_{\ell+p-1}}$.
By the above mentioned qc-lemma, $g_t$ is quasi-conformally conjugate to some rational function
 $$R_t=h_t\circ g_t\circ h_t^{-1}.$$
The quasi-conformal mapping $h_t:\widehat\C\longrightarrow\widehat\C$ is uniquely determined by the normalisation $h_t(z)=z$ for $z=0,\alpha,1$,  and depends
analytically on the parameter $t$. Also $h_t$ is analytic on $\widehat\C\setminus\overline{\bigcup_{p\ge 0}f_{t_0}^{-p}(A)}$, which, in particular, contains the points $0$, $v^\alpha_t$, and $v^\alpha_{t_0}$. We set $z_0=h_t(-1)$ to obtain $R_t(z)=a(t)z^m\big(\frac{1-z}{z-z_0}\big)^n.$
Since $h_t(\alpha)=\alpha$, $R_t$ has a critical point at
$z=\alpha$, and solving $R_t'(\alpha)=0$ for $z_0$ yields $z_0=-1$,
 thus
$$R_t(z)=a(t)z^m\Big(\frac{1-z}{1+z}\Big)^n.$$
From $R_t=h_t\circ\eta _t\circ f_{t_0}$ and $h_t(\alpha)=\alpha$,
however, it follows that
 $$a(t)v^\alpha_1=R_t(\alpha)=h_t\circ\eta_t\circ f_{t_0}(\alpha)=h_t\circ\eta_t(v^\alpha_{t_0})=h_t(v^\alpha_t),$$
hence $R_t(z)=f_\tau(z)$ with
$\tau=\tau(t)={h_t(v^\alpha_t)}/{v^\alpha_1}$ and $v^\alpha_\tau=h_t(v^\alpha_t);$
in particular, $\tau$ depends analytically on $t$. On some neighbourhood of $z=0$ we have
$$\begin{array}{rcl}
\ds(t_0/\tau)^{\frac{1}{m-1}}\Phi_{t_0}\circ h_t^{-1}\circ f_\tau &=&\ds(t_0/\tau)^{\frac{1}{m-1}}\Phi_{t_0}\circ
g_t\circ h_t^{-1}\cr
&=&\ds(t_0/\tau)^{\frac{1}{m-1}}\Phi_{t_0}\circ\eta_t\circ f_{t_0}\circ h_t^{-1}\cr
&=&\ds(t_0/\tau)^{\frac{1}{m-1}}\Phi_{t_0}\circ f_{t_0}\circ h_t^{-1}\cr
&=&\ds (t_0/\tau)^{\frac{1}{m-1}}t_0(\Phi_{t_0}\circ h_t^{-1})^{m}\cr&=&\ds \tau((t_0/\tau)^{\frac{1}{m-1}}\Phi_{t_0}\circ  h_t^{-1})^{m},
\end{array}$$
hence $\phi_\tau=(t_0/\tau)^{\frac{1}{m-1}}\Phi_{t_0}\circ h_t^{-1}$ solves B\"ottcher's functional equation
$$\phi_\tau\circ f_\tau(z)=\tau(\phi_\tau(z))^m.$$
Since $\tau$ and $h_t$ depend analytically on $t$, this is also true for $h_t^{-1}$, which is not self-evident.
Also from $h_t(g_t(z))=f_\tau(h_t(z))\sim\tau h_t(z)^m$ and $g_t(z)=f_{t_0}(z)\sim t_0z^m$ as $z\to 0$ it follows that $h_t(t_0z^m)\sim \tau h_t(z)^m$, hence
$h_t(z)\sim\sqrt[m-1]{t_0/\tau}z$, $h^{-1}_t(z)\sim\sqrt[m-1]{\tau/t_0}z$ and $\phi_\tau(z)\sim \lambda z$ as $z\to 0$, with $\lambda^{m-1}=1$. This implies $\phi_\tau=\lambda\Phi_\tau$
by uniqueness of the B\"ottcher coordinate, and from $\tau(t_0)=t_0$ and analytic dependence on $t$  it follows that $\lambda=1$ and $\phi_\tau=\Phi_\tau$, hence
$$\begin{array}{rcl}
E_k(\tau)&=&\tau^{\frac{1}{m-1}}\Phi_\tau(Q_k(\tau))=\tau^{\frac{1}{m-1}}\Phi_\tau(f^k_\tau(v^\alpha_\tau))\cr
&=& t_0^{\frac{1}{m-1}}\Phi_{t_0}\circ h_t^{-1}(f^k_\tau(v^\alpha_\tau))=t_0^{\frac{1}{m-1}}\Phi_{t_0}\circ
f^k_{t_0}\circ h_t^{-1}(v^\alpha_\tau)\cr
&=&t_0^{\frac{1}{m-1}}\Phi_{t_0}(f^k_{t_0}(v^\alpha_t))\quad{\rm if~} k\ge 1,~{\rm and}\cr
E_0(\tau)&=&t_0(\Phi_{t_0}(v^\alpha_t))^m.
\end{array}$$
Since $t\mapsto\tau$ is locally univalent, $E_k$ is univalent at $t_0$ if and only if the map $t\mapsto t_0^{\frac1{m-1}}\Phi_{t_0}(f_{t_0}^k(v^\alpha_t))$ is univalent on some neighbourhood of $t_0$. If $k\ge 1$, $\Phi_{t_0}$ is univalent on $\A_{t_0}$, and $f_{t_0}^k$ is univalent on $|z-v^\alpha_{t_0}|<\delta$ provided $Q_k(t_0)=f^k_{t_0}(v^\alpha_{t_0})\ne 0$, while $f_{t_0}^k$
is $n$-valent at $v^\alpha_{t_0}$ if $Q_k(t_0)=0$.
In case of $k=0$ we note that $\Phi_{t_0}$ is locally univalent on some forward  invariant domain $D$ that contains $0$ and $v^\alpha_{t_0}$,
and $v^\alpha_t=tv^\alpha_1\ne 0$ is trivially univalent.\hfill$\Box$

\medskip The proof of Theorem~\ref{SimpleConn} will be finished by

\begin{lem}\label{Lemma3}Let $h$ be a proper map of degree $m$ of the domain $D$ onto the unit disc $\D$, and assume that $h$ is ramified exactly over zero, that is, $h'(z)=0$ implies $h(z)=0$.
Then $D$ is simply connected and $h$ has a single zero on $D$.\end{lem}

\proof Assume that $h$ has zeros with multiplicities $m_\nu$ ($1\le\nu\le n)$. Then $h$ has degree $d=m_1+\cdots+m_n$ and $r=d-n$ critical points. The Riemann-Hurwitz formula
then yields $\# D-2=-d+r=-n$, hence $\#D=2-n$, which only is possible if $n=1$ and $\# D=1$.\hfill$\Box$

\medskip{\sc Remark.} Each connected component of $\OM^\alpha_k$ contains a zero of $Q_k(t)=tf_1(Q_{k-1}(t))$ which is not a zero of $Q_{k-1}$, hence is a zero of $Q_{k-1}(t)-1$.
Thus $\OM^\alpha_k$ consists of at most $\frac{d^k-1}{d-1}$ connected components.

\section{The hyperbolic loci}
\noindent The bifurcation locus $\BB^\beta$ is contained in some annulus $\delta<|t|<1$, and this also holds for $\W^\beta$. Hence (super-)attracting cycles $U_1,\ldots,U_k$
that contain the critical point $\beta$ may occur only if $\delta<|t|<1$.

\begin{thm}\label{polylike}For $0<|t|<1$,  $f_t$ is quasi-conformally conjugate to some polynomial
$$P_c(z)=cz^m(z+1)^n\quad(c=c_t\ne 0)$$
with free critical point $-\frac{m}{m+n}$. The basin $\A_t$ is completely invariant, and simply connected if and only if $t\not\in\OM^\alpha_{\rm res}.$
For $t\not\in\OM^\beta_0$, the Fatou set consists of $\A_t$, the simply connected basin $\B_t$ and
its pre-images and, additionally, of some (super-)attracting cycle of Fatou components and their pre-images if $t\in\W^\beta$; the cycle absorbs the critical point $\beta$.
\end{thm}

\proof To prove the second assertion we note that by Lemma~\ref{TRENNUNG} the pre-image $D$ of the disc $\triangle=\triangle_{r|t|}$ is a simply connected Jordan domain that contains $\overline\triangle\cup[0,1]$, but does not contain $v^\beta_t$. Then $D_2=\widehat\C\setminus\overline\triangle$ is a backward invariant domain, and
$$f_t:D_1\kzu1{d:1}D_2\quad(D_1=f_t^{-1}(D_2))$$
is a polynomial-like mapping in the sense of \cite{DouHub}, of degree $d=m+n$, hence is hybrid equivalent to some polynomial $P$ of degree $d$. We may assume that the quasi-conformal conjugation $\psi_t$ with
$$\psi_t\circ f_t=P\circ \psi_t$$
maps $\infty, 0,$ and $ -1$ onto $0,\infty,$ and $-1$, respectively. Thus $P$ is given by $P(z)=P_c(z)=cz^m(z+1)^n$, and $\psi_t$, hence also $c=c_t$ depends analytically on $t$.
\hfill$\Box$

\medskip{\sc Remark.} We note that $D_2=D_2(|t|)=\{z:|z|>r|t|\}$ increases if $|t|$ decreases, while $D_1=f_t^{-1}(\widehat\C\setminus\overline\triangle_{r|t|})=f_1^{-1}(\widehat\C\setminus\overline\triangle_r)$ is independent of $t$. Thus the conformal modulus $\mu(|t|)$ of $D_2(|t|)\setminus\overline{D_1}$
satisfies $\mu(1)\le\mu(|t|)-\log\frac1{|t|}\le\log\frac{\inf_{z\in D_1}|z|}{r}$.
The bifurcation locus of $P_c$ corresponds conformally to the bifurcation locus $\BB^\beta$, and the hyperbolic components are just quasi-conformal images of the hyperbolic components of the quadratic family $z\mapsto z^2+\xi$.

For $t\in\W_k$, the multiplier map $t\mapsto\lambda_t$ is an algebraic function of $t$. This is easily seen by writing the equations $f^k_t(z)=z$ and $\lambda=(f^k_t)'(z)$ as polynomial equations $q_1(z,t)=0$ and $q_2(z,t,\lambda)=0$, and computing the resultant $R_f(t,\lambda)$ of $q_1$ and $q_2$ with respect to $z$. For example, in case of $k=1$, $m=2$, and $n=1$ we obtain
$$R_f(t,\lambda)=[-2+14t-2t^2]+[1-10t+t^2]\lambda+2t\lambda^2=0.$$
For $P_c(z)=cz^2(z+1)$ we obtain in the same manner (multiplier $\mu$)
$$R_P(c,\mu)=9+2c-(c+6)\mu+\mu^2=0.$$
Since the quasi-conformal conjugation respects multipliers ($\lambda_t=\mu_{c_t})$, $c_t$ is an algebraic function of $t$ by the identity theorem; in the present case we obtain
$(1+2t+t^2+2tc)^2=0$ by computing the resultant of $R_f(t,\lambda)$ and $R_P(c,\lambda)$ with respect to $\lambda$, hence
$$t\mapsto c=c_t=-{\textstyle\frac12}\Big(t+2+\frac1t\Big)\quad(\textstyle c=-\frac 92\leftrightarrow t=\frac12(\sqrt{49}-\sqrt{45}))$$
maps $0<|t|<1$ conformally onto $\C\setminus[-2,0]$.

\begin{center}
\includegraphics[scale=0.17]{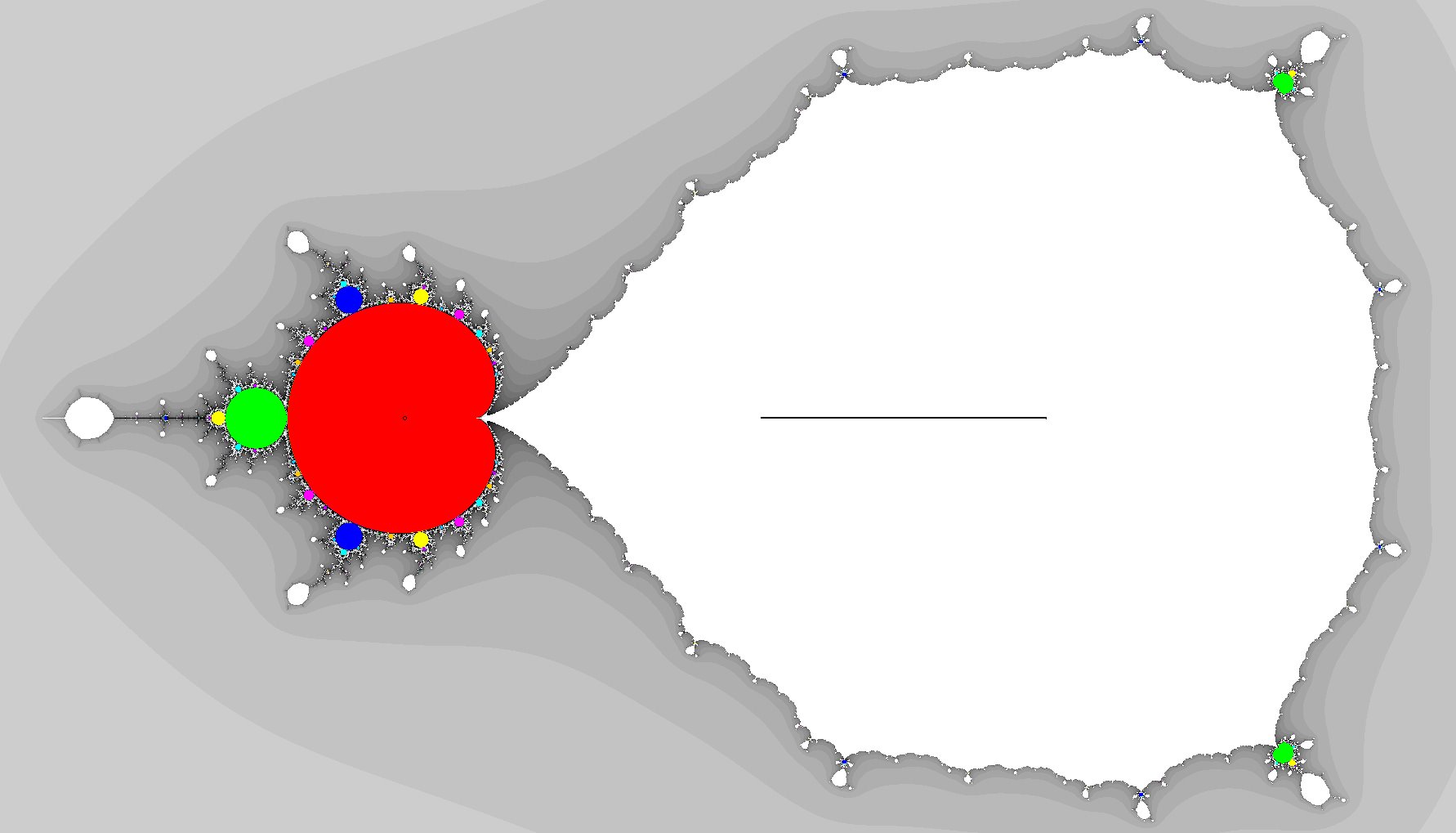}~\includegraphics[scale=0.17]{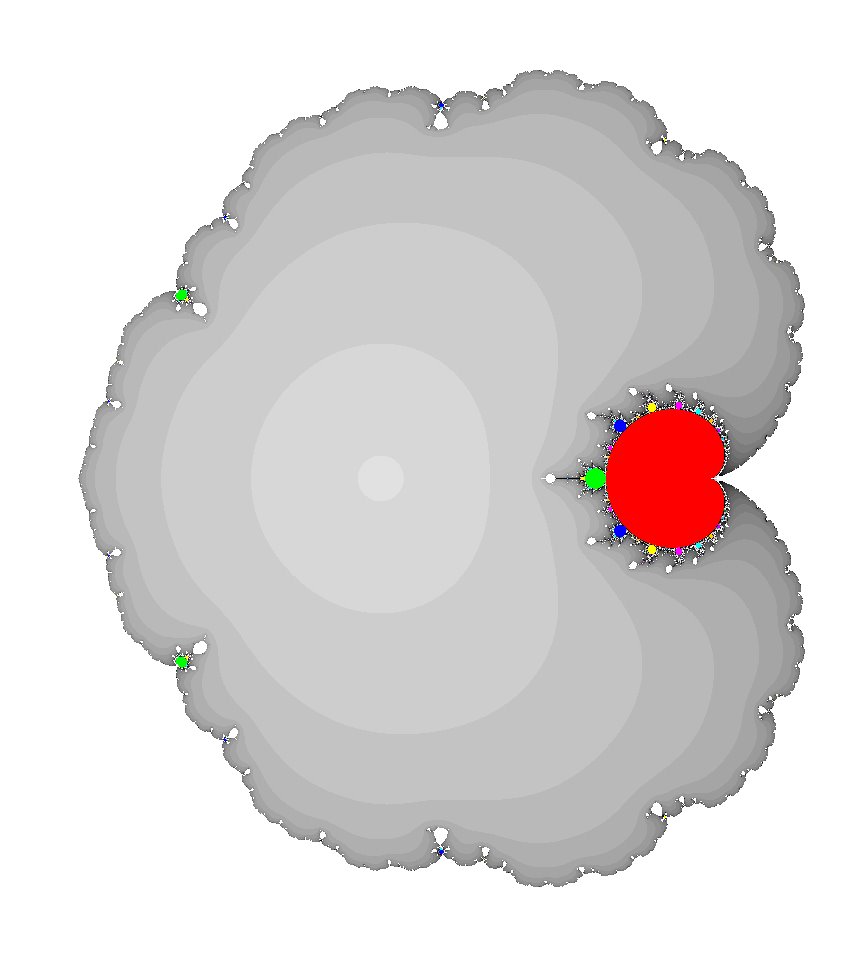}

\medskip\footnotesize{\sc Figure 2}-left: The parameter plane of $P_c(z)=cz^2(z+1)$. The escape region for $P_c$ (gray), the white region with slit, and the coloured domains correspond to $\OM^\alpha_{\rm res}$, $\OM^\beta\cap\D$, and $\W^\beta$, respectively. The punctured disc $0<|t|<1$ corresponds to $\C\setminus[-2,0]$ in the $c$-plane. {\sc Figure 2}{-right:} The parameter plane of $P_{-\frac12(t+2+\frac1t)}$ in $-0.2<\Re t<0.25,$ $-0.25<\Im t<0.25$ (see also Figure 1-right).
\end{center}

The following result was not explicitly stated but proved in \cite{Ste3}. The proof is an adaption of the procedure due to  Douady \cite{Douady}, 
applied to the hyperbolic components of the quadratic family $R_t(z)=z^2+t$ with one free critical point. The occurrence of several critical points requirers a slightly more sophisticated argument. The present version applies to a wider class of functions like $R_t(z)=z^d+t$, $R_t(z)=z^m+t/z^n$, $R_t(z)=t\big(1+\frac{((d-1)^{d-1}/d^d)z^d}{1-z}\big)$ ($d\ge 3$), $R_t(z)=-\frac t4\frac{(z^2-2)^2}{z^2-1}$, the present family, and many others, to show that the hyperbolic components are simply connected and are mapped properly onto the unit disc by the multiplier map $t\mapsto\lambda_t$.

\begin{thm}Let $(R_t)_{t\in T}$ be any family of rational maps that is analytically parametrised over some domain $T$. Suppose that each $R_t$ has a {(super-)}attracting cycle
$U_0\kzu1{m_1:1}U_1\kzu1{m_2:1}\cdots\kzu1{m_{n-1}:1}U_{n-1}\kzu1{m_{n}:1}U_{n}=U_0,$
such that $R_t^n$ has a single critical point $c_t\in U_0$ of multiplicity $m-1$, where $m=m_1\cdots m_n$ is the degree of $R^n_t:U_0\kzu1{m:1} U_0$. Assume also that the multiplier $\lambda_t$
satisfies $|\lambda_t|\to 1$ as $t\to\partial T$.
Then  the multiplier map $t\mapsto\lambda_t$ provides a proper map $T\kzu1{(m-1):1}\D$ which is ramified just over $w=0$, and $T$ is simply connected.\end{thm}

\end{document}